\documentclass[11pt]{article}
\font\tenmath=msbm10
\font\sevenmath=msbm7
\font\fivemath=msbm5
\newfam\mathfam \textfont\mathfam=\tenmath
\scriptfont\mathfam=\sevenmath \scriptscriptfont\mathfam=\fivemath

\usepackage{amsmath, moreverb, hyperref}
 \hypersetup{colorlinks=true,citecolor=blue}

\usepackage{natbib, color}
\usepackage{amsfonts, fancybox}
\usepackage{amssymb}
\usepackage[american]{babel}
\usepackage{graphicx}
\usepackage{flafter}
\usepackage[section]{placeins}

\newcommand{\cC}{\mathcal{C}}

\newcommand{\cF}{\mathcal{F}}

\newcommand{\cJ}{\mathcal{J}}
\newcommand{\cK}{\mathcal{K}}

\newcommand{\cQ}{\mathcal{Q}}

\newcommand{\cX}{\mathcal{X}}

\newcommand{\R}{{\rm I}\kern-0.18em{\rm R}}
\newcommand{\h}{{\rm I}\kern-0.18em{\rm H}}
\newcommand{\K}{{\rm I}\kern-0.18em{\rm K}}
\newcommand{\p}{{\rm I}\kern-0.18em{\rm P}}
\newcommand{\E}{{\rm I}\kern-0.18em{\rm E}}
\newcommand{\Z}{{\rm Z}\kern-0.18em{\rm Z}}
\newcommand{\1}{{\rm 1}\kern-0.24em{\rm I}}
\newcommand{\N}{{\rm I}\kern-0.18em{\rm N}}

\newcommand{\pn}{\p_{\kern-0.25em n}}
\newcommand{\pnm}{\p_{\kern-0.25em n,m}}
\newcommand{\psubm}{\p_{\kern-0.25em m}}
\newcommand{\psubp}{\p_{\kern-0.25em p}}
\newcommand{\cfi}{\cF_{\kern-0.25em \infty}}

\newcommand{\sgn}{\mathrm{sign}}

\newcommand{\argmin}{\mathop{\mathrm{argmin}}}
\newcommand{\argmax}{\mathop{\mathrm{arg\ max}}}
\newcommand{\ud}{\mathrm{d}}

\newcommand{\epr}{\hfill\hbox{\hskip 4pt\vrule width 5pt
                  height 6pt depth 1.5pt}\vspace{0.5cm}\par}

\newcommand{\bi}{^{(i)}}
\newcommand{\hpi}{\hat \pi}
\newcommand{\spi}{\pi^{\star}}

\newcommand{\stn}{\sum_{t=1}^n}

\newcommand{\sst}{\sum_{s=1}^t}
\newcommand{\sjM}{\sum_{j=1}^{M^d}}

\newtheorem{TH1}{Theorem}[section]

\newtheorem{prop}{Proposition}[section]
\newtheorem{cor}{Corollary}[section]

\newtheorem{lem}{Lemma}[section]

\begin{document}

\title{Nonparametric Bandits with Covariates}

\author{
{\sc Philippe Rigollet}
 \thanks{Princeton University. Partially supported by the National Science Foundation (DMS-0906424).} 
 \quad  \quad  {\sc Assaf Zeevi} \thanks{Columbia University.}
}

\date{\normalsize \today}

\maketitle

\begin{abstract}
We consider a bandit problem which involves sequential sampling 
from two populations (arms). Each arm produces a noisy reward
realization which depends on an observable random {\it covariate}.
The goal is to maximize cumulative expected reward. We derive
general lower bounds on the performance of any admissible policy,
and develop an algorithm
 whose performance achieves the order of said lower bound up to logarithmic terms. This is done by decomposing the global problem into
 suitably ``localized'' bandit
 problems.
 Proofs blend ideas from nonparametric statistics and traditional
 methods used in the bandit literature.
\end{abstract}
\medskip



\section{Introduction}
The seminal paper of \citet{Rob52} introduced an important
class of sequential optimization problems, otherwise known as
multi--armed bandits. These models have since been used
extensively in such fields as statistics, operations research,
engineering, computer science and economics. The traditional
two--armed bandit problem can be described as follows. Consider
two statistical populations (arms), where at each point in time it
is possible to sample from only one of the two and receive a
random reward dictated by the properties of the sampled
population. The objective is to
 devise a sampling policy that maximizes  expected cumulative
 (or discounted) rewards over a finite (or infinite)
 time horizon.  The difference between the performance of said
 sampling policy and that of an oracle, that repeatedly samples
 from the population with the higher mean reward, is called the
 {\it regret}. Thus, one can re-phrase the objective as minimizing
 the regret.
\par
The original motivation for bandit-type problems originates from
treatment allocation in clinical trials; see, e.g., \citet{LaiRob85}
for further discussion and references therein.  Here patients
enter sequentially and receive one of several treatments. The
efficacy of each treatment is unknown, and for each patient a
noisy measurement of it is recorded. The goal is to assign as many
patients as possible to the best treatment. An example of more
recent  work can be found in the area of web-based advertising,
and more generally
 customized marketing. An on-line
  publisher needs to choose one of
 several ads to present to consumers, where the
  efficacy of these ads is unknown.  The publisher observes
  click-through-rates (CTRs) for each ad, which provide a noisy measurement of the
  efficacy, and based on that needs to assign ads that
  maximize CTR.
\par
When the populations being sampled are homogenous, i.e., when the
sequential rewards are independent and identically distributed
(iid) in each arm, \citet{LaiRob85} proposed  a family of
policies that at each step compute the empirical mean reward in
each arm, and adds to that a confidence bound that accounts for
uncertainty in these estimates. These so-called
upper-confidence-bound (UCB) policies were shown to be
asymptotically optimal. In particular, it is proven in
\citet{LaiRob85} that such a policy incurs a regret of order
$\log n$, where $n$ is the length of the time horizon, and no
other ``good" policy can (asymptotically) achieve a smaller
regret; see also \citet{AueCesFis02}. The elegance of the theory
and sharp results developed in \citet{LaiRob85} hinge to a large
extent on the assumption of homogenous populations and hence
identically distributed rewards. This, however,  is clearly too
restrictive for many applications of interest. Often, the decision
maker observes further information and based on that a more {\it
customized} allocation can be  made. In such settings rewards may
still be assumed to be
 independent,  but no longer identically
distributed in each arm. A particular way to encode this is to
allow for an exogenous variable (a covariate) that affects the
rewards generated by each arm at each point in time when this arm
is pulled.
\par
Such a formulation was first introduced in
\citet{Woo79} under parametric assumptions and in a somewhat
restricted setting; see \citet{GolZee09} and \citet{WanKulPoo05} for  two
  very  different recent approaches to the study of
  such bandit problems, as well as references therein for further
  links to antecedent literature. The first work to venture outside
  the realm of parametric modeling assumptions was that of
  \citet{YanZhu02}. In particular, they assumed the mean response in each arm,
  conditional on the covariate value, follows a general
  functional form, hence one can view their setting as as {\it
  nonparametric}
  bandit problem.  They proposed a policy that is based on
  estimating each response function, and then, rather than
 greedily  choosing the arm with the highest estimated mean
 response given the covariate,
  allows with some small probability of selecting a potentially
  inferior arm. (This is a variant of $\varepsilon$-greedy
  policies; see \citet{AueCesFis02}.)
  If the nonparametric estimators of the arms' functional response
  are consistent, and the randomization is chosen in a suitable manner,
    then the above policies ensure that the average
  regret tends to zero as the time horizon $n$ grows to infinity.
  In the typical bandit terminology, such policies are said to be
  {\it consistent}.
 However, it is unclear whether they
   satisfy a more refined notion of optimality, insofar as the magnitude of the regret
   is concerned, as is the case for
   UCB-type policies in traditional bandit problems. Moreover,
   the study by \citet{YanZhu02} does not spell out the
   connection between the characteristics of the class of response
   functions, and the resulting complexity of the nonparametric
   bandit problem.
     \par
  The purpose of the present paper is to further understanding of
  nonparametric bandit problems, deriving regret-optimal
 policies and shedding light on some of the elements that
 dictate the complexity of such problems.
   We make only two
assumptions on the underlying functional form that governs the
arms' responses. The first is a mild smoothness condition.
Smoothness assumptions can be exploited using  ``plug-in" policies
as opposed ``minimum contrast" policies; a detailed account of the
differences and similarities between these two setups in the full
information case can be found in~\citet{AudTsy07}.   Minimum
contrast type policies have already received some attention in the
bandit literature with side information, aka \emph{contextual
bandits}, in the papers of \citet{LanZha08} and also
\citet{KakShaTew08}. In these studies, admissible policies are
restricted to a more limited set than the general class of
non-anticipating policies.  A related problem online convex
optimization with side information was studied by
\citet{HazMeg07}, where the authors use discretization technique
similar to the one employed in this paper. It isi worth noting
that the cumulative regret in these papers  is defined in a weaker
form compared to the traditional bandit literature, since the
cumulative reward of a proposed policy is compared to that of the
best policy in a certain restricted class of policies. Therefore,
 bounds on the regret depend, among other things, on the complexity of said class of
policies. Plug-in type policies have received attention in the
context of the continuum armed bandit problem, where as the nsame
suggests  there are uncountably many arms. Notable entries in that
stream of work are \citet{Sli09} and \citet{LuPalPal09}, who
impose a smoothness condition  both on the space of arms and the
space of covariates, obtaining optimal regret bounds up to
logarithmic terms.
\par
The second key assumption  in our paper is a so-called  {\it
margin condition}, as it has been come be known  in the full
information setup; cf. \citet{Tsy04}.  In that setting, it has
been shown to critically affect  the complexity classification
problems ~\cite{Tsy04,BouBouLug05,AudTsy07}. In the bandit setup,
this condition encodes the ``separation'' between the functions
that describe the arms' responses and was originally studied
by~\citet{GolZee09} in the one armed bandit problem; see further
discussion in section \ref{sec:desc}. We will see later that the
margin condition is a natural measure of complexity in the
nonparametric bandit problem.
\par
In this paper, we introduce a family of policies called UCBograms.
The term is indicative of two salient ingredients of said
policies: they build on {\it regressogram} estimators;  and
augment the resulting mean response estimates with
upper-confidence-bound terms. The idea of the regressogram is
quite natural and easy to implement. It groups the covariate
vectors into bins and then estimates, by means of
 simple averaging, a constant which is a proxy for the mean
 response of each arm over each such bin.
One then views these bins as indexing ``local'' bandit problems,
which are solved by applying a suitable UCB-type modification,
following the logic of \citet{LaiRob85} and \citet{AueCesFis02}.
In other words, this family of policies decomposes the
non-parametric bandit problem into a sequence of localized
standard bandit problems; see section \ref{SEC:main} for a
complete description. The idea of binning covariates lends itself
to natural implementation in the two motivating examples described
earlier: patients and consumers are segmented into groups with
``similar'' characteristics; and then the treatment or ad is
allocated based on the characteristic response over that group.
\par
In terms of performance, we  prove that the UCBogram policies
achieve a regret that is fairly large compared to typical orders
of regret observed in the literature. In particular, as opposed to
a bounded or logarithmic growth, in our setting the order of the
regret is {\it polynomial} in the time horizon $n$; see Theorem
\ref{TH:main}. One may question, especially given the simple
structure and logic underlying the UCBogram policy, whether this
is the best that can be achieved in such problems. To that end, we
prove a lower bound which demonstrates that for any admissible
policy there exist arm response functions satisfying our
assumptions for which one cannot improve on the polynomial order of the upper
bound established in Theorem \ref{TH:main};  see Theorem
\ref{th:minimax}. Finally, beyond these analytical results, in our
view one of the contributions of the present paper is in pointing
to some possible synergies and potentially interesting connections
between  the traditional bandit literature and nonparametric
statistics.

\section{Description of the problem}
\label{sec:desc}

\subsection{Machine and game}
A \emph{bandit machine with covariates} is
 characterized by a sequence
$$
(X_t, Y_t^{(1)}, Y_t^{(2)}), \ t=1,2,\dots
$$
of independent random vectors, where $\big(X_t\big)$,
$t=1,2,\ldots$ is a sequence of iid covariates in $\cX \subset
\R^d$ with probability distribution $P_X$,  and $Y_t\bi$ denotes
the random reward yielded by arm $i$ at time $t$. We assume that,
for each $i=1,2$, conditionally on $\{X_t=j\}$, the rewards
$Y_t\bi, t=1, \ldots, n$ are i.i.d random variables in $[0,1]$
with conditional expectation given by
$$
\E\big[Y_t\bi|X_t] = f\bi(X_t)\,, \quad t=1,2,\dots, \ i=1,2\,,
$$
where $f\bi, i=1,2$, are unknown functions such
that $0 \le f\bi(x) \le 1$, for any $i=1,2,\  x \in \cX$. A natural example arises when $Y_t\bi$ takes values in $\{0,1\}$ so that the conditional distribution of $Y_t\bi$ given $X_t$ is Bernoulli with parameter $f\bi(X_t)$.

The \emph{game} takes place sequentially on this machine, pulling
one of the two arms at each time $t=1, \ldots, n$. A
 {\it non-anticipating policy} $\pi=\{\pi_t\}$ is a sequence
of random functions $\pi_t:\cX \to\{1,2\} $ indicating to the
operator which arm to pull at each time $t$, and such that $\pi_t$
depends only on observations strictly anterior to~$t$. The
\emph{oracle rule} $\pi^\star$, refers to the strategy that would
be played by an omniscient operator with complete knowledge of the
functions $f\bi, i=1,2$. Given side information $X_t$, the oracle
policy $\spi$ prescribes the arm with the largest expected reward,
i.e.,
$$
\pi^\star(X_t):=\argmax_{i=1,2} f\bi(X_t)\,.
$$
The oracle rule will be used to benchmark any proposed policy
$\pi$ and to measure the performance of the latter via its
\emph{(expected cumulative) regret} at time $n$ defined by
$$
R_n(\pi):=\E \sum_{t=1}^n \big(Y^{(\pi^\star(X_t))}_t-Y^{(\pi_t(X_t))}_t \big)=\E \sum_{t=1}^n \big(f^{(\pi^\star(X_t))}(X_t)-f^{(\pi_t(X_t))}(X_t)\big)\,.
$$
Also, let $S_n(\pi)$ denote the {\it inferior sampling rate} at time
$n$ defined by
\begin{equation}\label{e-isr}
S_n(\pi):=\E\stn \1(\pi_t(X_t)\neq \spi_t(X_t), f^{(1)}(X_t)\neq
f^{(2)}(X_t))\,,
\end{equation}
where $\1(A)$ is the indicator function that takes value $1$ if
event $A$ is realized and $0$ otherwise. The quantity $S_n(\pi)$
measures the expected number of times at which a \emph{strictly}
suboptimal arm has been pulled, and note that in our setting the
suboptimal arm varies as a function of the covariate value $x$.

Without further assumptions on the machine, the game can be
arbitrarily difficult and, as a result, the regret and inferior
sampling rate can be arbitrarily close to $n$. In the following
subsection, we describe natural  assumptions on the regularity of
the machine that allow to control its complexity.

\subsection{Smoothness and margin conditions}

As usual in nonparametric estimation we first impose some
regularity on the functions $f\bi, i=1,2$. Here and in what
follows we use $\|\cdot\|$ to denote the Euclidean norm.

\medskip

\noindent {\sc Smoothness condition.} We say that the machine
satisfies the smoothness condition with parameters $(\beta, L)$ if
\begin{equation}
\label{EQ:holder}
|f\bi(x)-f\bi(x')|\le L\|x-x'\|^\beta, \quad \forall\,  x, x' \in \cX, i=1,2
\end{equation}
for some $\beta \in (0,1]$ and $L>0$.

\medskip

Notice that a direct consequence of the smoothness condition with
parameters $(\beta, L)$ is that the function
$\Delta:=|f^{(1)}-f^{(2)}|$ also satisfies the smoothness
condition with parameters $(\beta, 2L)$. The behavior of function
$\Delta$ critically controls the complexity of the problem and the
smoothness condition gives a local upper bound on this function.
 The second condition
 imposed gives a  lower bound on this function though in a weaker
global sense. It is closely related to the margin condition
employed in classification \cite{Tsy04,MamTsy99}, which drives the
terminology employed here.

\medskip

\noindent {\sc Margin condition.}  We say that the machine
satisfies the margin condition with parameter $\alpha$ if there
exists $\delta_0\in(0,1)$, $C_\delta>0$ such that
$$
P_X\big[ \,0< |f^{(1)}(X)-f^{(2)}(X) |\le \delta\big]\le C_\delta
\delta^\alpha\,, \quad \forall\,   \delta \in [0, \delta_0]\,
$$
 for some  $\alpha>0$.
\medskip

In what follows, we will focus our attention on marginals $P_X$
that are equivalent to the Lebesgue measure on a compact subset of
$\R^d$. In that way, the margin condition will only contain
information about the behavior of the function $\Delta$ and not
the marginal $P_X$ itself. A large value of the parameter $\alpha$
means that the function $\Delta$ either takes value 0 or is
bounded away from 0,  except over a set of small
$P_X$-probability. Conversely, for values of $\alpha$ close to 0,
the margin condition is essentially void and the two functions can
be arbitrary close, making it difficulty to distinguish among
them. This will be reflected in the bounds on the regret which are
derived in the subsequent section.

Intuitively, the smoothness condition and the margin condition
work in opposite directions.  Indeed, the former ensures that the
function $\Delta$ does not depart from zero too fast whereas the
latter warrants the opposite. The following proposition accurately
quantifies the extent to which the conditions are conflicting.

\begin{prop}
\label{prop:alpha-beta} Under the smoothness condition with
parameters $(\beta, L)$, any machine that satisfies the margin
condition with parameter $\alpha$ such that $\alpha \beta>1$
exhibits an oracle policy $\spi$ which dictates pulling only one
of the two arms all the time,  $P_X$-almost surely. Conversely, if
$\alpha\beta \le 1$ there exists machines with nontrivial oracle
policies.
\end{prop}
{\bf Proof.} The first part of the proof is a straightforward
consequence of Proposition~3.4 in \citet{AudTsy07}. To prove the
second part, consider the following example. Assume that $d=1$,
$\cX=[0,2]$, $f^{(2)}\equiv 0$ and  $f^{(1)}(x)=L{\rm
sign}(x-1)|x-1|^{1/\alpha}$. Notice that $f^{(1)}$ satisfies the
smoothness condition with parameters $(\beta, L)$ if and only if
$\alpha\beta \le 1$. The oracle policy is not trivial and defined
by $\spi(x)=2$ if $x\le 1$ and $\spi(x)=1$ if $x> 1$. Moreover, it
can be easily shown that the machine satisfies the margin
condition with parameter $\alpha$ and with  $\delta_0=C_\delta=1$.
\epr

\section{Policy and main result}

\label{SEC:main} We first outline a policy to operate the
bandit machine described in the previous section. Then we state
the main result which is an upper bound on the regret for this
policy. Finally, we state a proposition which allows us to
translate the bound on the regret into a bound on the inferior
sampling rate.

\subsection{Binning and regressograms}

To design a policy that solves the bandit problem described in the
previous section, one has to inevitably find an estimate of the
functions $f\bi, i=1,2$ at the current point $X_t$. There exists a
wide variety of nonparametric regression estimators ranging from
local polynomials to wavelet estimators. However, a very simple
piecewise constant estimator, commonly referred to as
\emph{regressogram} will be particularly suitable for our
purposes.

Assume now that $\cX=[0,1]^d$ and let $\{B_j, j=1,\dots, M^d\}$ be the
regular partition of $\cX$, i.e., the reindexed collection of
hypercubes defined for ${\sf k}=(k_1, \ldots, k_d)\in \{1, \ldots, M\}^d$\,,
$$
B_{\sf k}=\Big\{x\in \cX\,:\, \frac{k_\ell-1}{M}\le
x_\ell\le \frac{k_\ell}{M}\,, \ell=1, \ldots, d\ \Big\}\,.
$$
For each arm $i=1,2$, consider the average reward for each bin
$B_j, j=1, \ldots, M^d$ defined by
$$
\bar f\bi_j=\frac{1}{p_j}\int_{B_j}f\bi(x)\ud x\,,
$$
where $p_j=P_X(B_j)$\,. By analogy with histograms, the empirical
counterpart of the piecewise constant function $x\mapsto \sjM \bar
f\bi_j\1(x\in B_j)$, is
often called \emph{regressogram}. To define it, we need the following
quantities. Let $N_t\bi(j,\pi)$ denote the
number of times $\pi$ prescribed to   pull arm $i$ at times anterior
to $t$ when the covariate
was in bin $B_j$,
$$
N_t\bi(j, \pi)=\sst \1(X_s\in B_j, \pi_s(X_s)=i)\,,
$$
and let $\overline  Y_t\bi (j, \pi)$ denote the average reward collected at those times,
$$
\overline  Y_t\bi (j, \pi)=\frac{1}{N_t\bi(j, \pi)}\sst Y\bi_s \1(X_s\in B_j,\pi_s(X_s)=i)\,,
$$
where here and throughout this paper, we use the convention  $1/0=\infty$.
For any arm $i=1,2$ and any time $t\ge 1$ the regressograms obtained
from a policy $\pi$ at time $t$ are defined by the following piecewise
constant estimators
$$
\hat f\bi_{t, \pi}(x)= \sjM \overline  Y_t\bi (j, \pi)\1(x\in B_j)\,.
$$
While regressograms are rather rudimentary  nonparametric
estimators of the functions $f\bi$, they allow us to decompose the
original problem into a collection of $M^d$ traditional bandit
machines without covariates, each one corresponding to a different
bin.

\subsection{The UCBogram}
\label{SUB:UCBogram}
The ``UCBogram'' is an index type policy based on upper
confidence bounds for the regressogram defined above. Upper confidence
bounds (UCB) policies are known to perform optimally in the
traditional two armed bandit problem, i.e., without covariates
\cite{LaiRob85,AueCesFis02}. The index of each arm is computed as the sum of the average past
reward and a stochastic term accounting for the deviations of the
observed average reward from the true average reward. In the UCBogram,
the average reward is simply replaced by the value of the regressogram
at the current covariate $X_t$.

For any $ s \ge 1$ the upper confidence bound at time $t$ bound is of the form
$$
U_t(s)=\sqrt{\frac{2\log t}{s}}\,.
$$
The UCBogram $\hpi$ is defined as follows. For any $x \in [0,1]^d$, define
$$
N\bi_t(x)=\sjM N\bi_t(j, \hpi)\1(x\in B_j)\,,
$$
the number of times the UCBogram prescribed to pull arm $i$ at times
anterior to $t$ when the covariate was in the same bin as $x$. Then $\hpi=(\hpi_1, \hpi_2, \ldots)$  is defined recursively by
$$
\hpi_t(x)=\argmax_{i=1,2}\Big \{\hat f\bi_{t, \hat \pi}(x) + U_{t}(N\bi_t(x))\Big\}\,.
$$
Notice that the UCBogram is indeed a UCB-type policy. Indeed, for
each arm $i=1, 2$ and at each point $x$, it computes an estimator
$\hat f\bi_{t, \pi}(x)$ of the expected reward and adds an upper
confidence bound $U_{t}(N\bi_t(x))$ to account for stochastic
variability in  this estimator. The most attractive feature of the
regressogram is that it allows to decompose the nonparametric
bandit  problem into independently operated local machines as
detailed in the proof of the following theorem.

\begin{TH1}
\label{TH:main} Fix $\beta \in (0,1]$, $L>0$ and $\alpha \in
(0,1]$. Let $\cX=[0,1]^d$ and assume that the covariates $X_t$
have a distribution which is
equivalent\footnote{\label{foot:equiv}Two measures $\mu$ and $\nu$
are said to be \emph{equivalent} if there exist two positive
constants $\underline{c}$ and $\bar c$ such that
$\underline{c}\mu(A) \le \nu(A) \le \bar c \mu(A)$ for any
measurable set $A$.} to the Lebesgue measure  on the unit
hypercube $\cX$. Let the machine satisfy both the smoothness
condition with parameter $(\beta, L)$ and the margin condition
with parameter $0<\alpha\le 1$. Then the UCBogram policy $\hpi$
with $M=\lfloor(n/\log n)^{1/(2\beta+d)}\rfloor$ has an expected
cumulative regret at time $n$ bounded by
$$
R_n(\hat \pi)\le Cn\max\Big\{\Big(\frac{n}{\log
  n}\Big)^{-\frac{\beta(\alpha+1)}{2\beta +d}}\,,\Big(\frac{n}{(\log
  n)^2}\Big)^{-\frac{2\beta}{2\beta+d}} \Big\}\,,
$$
where $C>0$ is a positive constant.
\end{TH1}
{\bf Proof.}
To keep track of positive constants, we number them $c_1, c_2, \ldots$. Define $c_1=2Ld^{\beta/2} +1,$ and let $n_0 \ge 2$ be the largest integer such that
$$
\left(\frac{n_0}{\log n_0}\right)^{\beta/(2\beta+d)} \le \frac{2c_1}{\delta_0}\,,
$$
where $\delta_0$ is the constant appearing in the margin condition.
If $n \le n_0$, we have $R_{n} \le n_0$ so that the result of the theorem holds when $C$ is chosen large enough, depending on the constant $n_0$.
In the rest of the proof, we assume that $n > n_0$ so that $c_1M^{-\beta} < \delta_0$.

Recall that the UCBogram
policy $\hpi$ is a collection of functions $\hpi_t$ that are constant on each
$B_j$, equal to $\hpi_t(j)$. Define the regret $\textsf{R}_j(\hpi)$ on bin $B_j$ by
$$
\textsf{R}_j(\hpi)= \stn \big(f^{(\pi^\star(X_t))}(X_t)-f^{(\hpi_t(j))}(X_t)\big)\1(X_t
\in B_j)\,,
$$
and observe that the overall regret of $\hpi$ can be written as
$$
R_n(\hpi)=
\sjM\E \textsf{R}_j(\hpi)\,.
$$
Consider the set of ``well behaved'' bins on which the expected reward
functions of the two arms are well separated:
$$
\cJ=\{j\,:\, \exists\ x \in B_j\,,\,
|f^{(1)}(x)-f^{(2)}(x)|> c_1M^{-\beta}  \}\,.
$$
For any  $j\notin \cJ$ and any $x \in B_j$, we have $|f^{(1)}(x)-f^{(2)}(x)|\le c_1M^{-\beta} < \delta_0$ so that
$$
\E \textsf{R}_j(\hpi)\le c_1M^{-\beta}\stn \p\big[0<|f^{(1)}(X_t)-f^{(2)}(X_t)|\le
c_1 M^{-\beta} , X_t
\in B_j\big]\,,
$$
Summing over $j \notin \cJ$, we obtain from the margin condition that
\begin{equation}
\label{eq:difficult}
\sum_{j\notin \cJ}\E \textsf{R}_j(\hpi) \le C_\delta c_1^{1+\alpha}n M^{-\beta(1+\alpha)}  \,.
\end{equation}

We now treat the well behaved bins, i.e., bins $B_j$ such that $j \in \cJ$. Notice that since each bin is a hypercube with side length $1/M$ and since the reward functions satisfy the smoothness condition with parameters $(\beta, L)$,
we have
$$
|f^{(1)}(x)-f^{(2)}(x)|> c_1M^{-\beta}-2Ld^{\beta/2}M^{-\beta}=M^{-\beta}\,,
$$
for
\emph{any} $x \in B_j, j \in \cJ$. In particular, for such $j$, since the two functions are
continuous, the difference $f^{(1)}(x)-f^{(2)}(x)$ has constant sign
over $B_j$ and $|\bar f^{(1)}_j-\bar f^{(2)}_j| >
M^{-\beta}$. As a consequence, the oracle policy $\spi$ is constant on
$B_j$, equal to $\spi(j)$ for any $j \in \cJ$ and, conditionally on $\{X_t\in B_j\}$, the game can be viewed as a standard bandit problem, i.e., without covariates, where arm $i$ has bounded reward with mean $\bar f\bi_j$. Moreover, conditionally on $\{X_t \in B_j\}$, the UCBogram can be seen as a standard UCB policy. Applying for example Theorem~1 in~\citet{AueCesFis02}, we find that for $j\in \cJ$,
\begin{equation}
\label{eq:easy}
\E \textsf{R}_j(\hpi)\le\Big[\big(1+\frac{\pi^2}{3}\big) \Delta_j\Big] +\frac{8\log n}{\Delta_j}\le c_2\frac{\log n}{\Delta_j}\,,
\end{equation}
where $\Delta_j=|\bar f^{(1)}_j-\bar f^{(2)}_j|$ is the average gap in bin~$B_j$. We now use the margin condition to provide lower bounds on $\Delta_j$.
Assume without loss of generality that the gaps are ordered
$0<\Delta_1 \le \Delta_2 \le \ldots, \le \Delta_{M^d}$ and define the integers $j_1, j_2$ such that $\cJ=\{j_1, \ldots, M^d\}$ and $j_2\in \{j_1, \ldots, M^d\}$ is the largest integer such that $\Delta_{j_2} \le \delta_0/c_1$. Therefore, for any $j \in \{j_1, \ldots, j_2\} \subset \cJ$, we have on the one hand,
\begin{equation}
\label{EQ:prUB1}
P_X\big[ \,0< |f^{(1)}-f^{(2)} |\le \Delta_j+(c_1-1)M^{-\beta}\big]\ge\sum_{k=1}^{M^d} p_k \1(0<\Delta_k \le \Delta_{j})\ge \frac{\underline{c} j}{M^d}\,,
\end{equation}
where we use the fact that $p_k=P_X(B_k) \ge \underline{c}/M^d$ since $P_X$ is equivalent to the Lebesgue measure on $[0,1]^d$ (see footnote~\ref{foot:equiv}). On the other hand, the margin condition yields for any $ j \in \{j_1, \ldots, j_2\}$ that,
\begin{equation}
\label{EQ:prUB2}
P_X\big[ \,0< |f^{(1)}-f^{(2)} |\le \Delta_j+(c_1-1)M^{-\beta}\big]\le C_\delta\big(c_1\Delta_j)^\alpha\,.
\end{equation}
where we used the fact that $\Delta_j+(c_1-1)M^{-\beta}\le c_1\Delta_j \le \delta_0$, for any  $j \in \{j_1, \ldots, j_2\}$. The previous two inequalities yield
\begin{equation}
\label{EQ:lowerDelta}
\Delta_{j}\ge c_3 \Big(\frac{j}{M^d} \Big)^{1/\alpha}\,, \quad \forall\ j \in \{j_1, \ldots, j_2\}\,.
\end{equation}
Combining (\ref{eq:difficult}), (\ref{eq:easy}) and (\ref{EQ:lowerDelta}), we obtain
the following bound,
\begin{equation}
\label{EQ:reg2}
R_n(\hpi)\le c_4\Big[nM^{-\beta(1+\alpha)}+ j_1M^{-\beta}+ (\log n)\sum_{j=j_1}^{j_2} \left(\frac{M^d}{j}\right)^{1/\alpha}  + M^d \log n\Big]\,.
\end{equation}
Note that applying the same arguments as in~\eqref{EQ:prUB1}
and~\eqref{EQ:prUB2}, we find that $j_1$ satisfies
$$
\frac{\underline{c} j_1}{M^d}\le P_X\big[ \,0< |f^{(1)}-f^{(2)} |\le c_1M^{-\beta}\big]\le C_\delta\big(c_1M^{-\beta})^\alpha\,,
$$
so that $j_1 \le c_5M^{d-\alpha\beta}$. We now bound from above the sum in ~(\ref{EQ:reg2}) using the following integral approximation:
\begin{equation}
\label{EQ:integral}
\sum_{j=j_1}^{j_2} \left(\frac{M^d}{j}\right)^{1/\alpha} \le \sum_{j=j_1}^{M^d}\left(\frac{M^d}{j}\right)^{1/\alpha} \le c_7 M^d \int_{M^{-\alpha\beta}}^{1}
 x^{-1/\alpha} \ud x\,.
\end{equation}
If $\alpha<  1$, this integral is bounded by $c_6M^{\beta(1-\alpha)} $ and if $\alpha =1$, it is
bounded by $c_7\log M$. As a result, the integral in~(\ref{EQ:integral}) is
of order $M^d(M^{\beta(1-\alpha)}\vee\log M)$ and we obtain from~(\ref{EQ:reg2}) that
\begin{equation}
\label{EQ:finalprUB}
R_n(\hpi)\le c_{8}\Big[nM^{-\beta(1+\alpha)}+ M^d(M^{\beta(1-\alpha)}\vee\log M)\log n
\Big]\,,
\end{equation}
and the result follows by choosing $M$ as prescribed. \epr

We should point out that the version of the UCBogram described
above specifies the number of bins $M$ as a function of the
horizon $n$, while in practice one does not have foreknowledge of
this value. This limitation can be easily circumvented by using
the so-called \emph{doubling argument} \cite{CesLug06} which
consists of ``reseting'' the game at times $2^{k}, k=1, 2, \ldots$

The reader will note that when $\alpha =1$ there is an additional
$\log n$ factor appearing in the upper bound given in the statement
of the theorem. More generally, for any $\alpha>1$, it is possible
to minimize the expression on the right hand side
of~\eqref{EQ:finalprUB} with respect to $M$, but the optimal value
of  $M$ would then depend on the value of $\alpha$. This sheds
some  light on a significant limitation of the UCBogram which
surfaces in this parameter regime: it requires the operator to
pull each arm at least once in each bin and therefore to incur a
regret of at least order $M^d$. In other words, the UCBogram
splits the space $\cX$ in ``too many'' bins when $\alpha\ge 1$.
Intuitively this can be understood as follows. When $\alpha=1$,
the gap function $\Delta(x)$ is bounded  away from zero for most
$x \in \cX$. For such $x$, there is no need to carefully estimate
the gap function since it has constant sign for ``large''
contiguous regions. As a result one could use larger bins in such regions reducing the overall number of bins and therefore removing the extra logarithmic term. Of course, such limitations are intrinsic to the UCBogram and may not appear with other policies but it is beyond the scope of this paper.

\subsection{The inferior sampling rate}

 Unlike traditional
bandit problems, the connection between the
inferior sampling rate defined in (\ref{e-isr}) and the regret is more intricate here.
The following lemma establishes a connections between the two.

\begin{lem}
\label{LEM:ISR}
For any $\alpha>0$, under the margin condition we have
$$
S_n(\pi)\le C
n^{\frac{1}{1+\alpha}}R_n(\pi)^{\frac{\alpha}{1+\alpha}}\,,
$$
for any policy $\pi$ and for some positive constant $C>0$.
\end{lem}
{\bf Proof.} The idea of the proof is quite standard and originally appeared in \citet{Tsy04}. It has been used in
\citet{RigVer09} and \citet{GolZee09}.  Define the two random quantities:
$$
r_n(\pi)=\stn |f^{(1)}(X_t)- f^{(2)}(X_t)|\1(\pi_t(X_t)\neq
\spi(X_t))\,,
$$
and
$$
s_n(\pi)=\stn \1(f^{(1)}(X_t)\neq f^{(2)}(X_t), \pi_t\neq \spi(X_t))\,.
$$
We have
\begin{eqnarray}
\label{EQ:ISR_pr1}
 r_n(\pi)&\ge& \delta \stn \1( \pi_t(X_t)\neq \spi(X_t))\1(|
 f^{(1)}(X_t)- f^{(2)}(X_t)|>\delta)\nonumber\\
& \ge &\delta\big[s_n(\pi)-  \stn \1(  \pi_t(X_t)\neq \spi(X_t),0<|
 f^{(1)}(X_t)- f^{(2)}(X_t)|\le \delta)\big]\nonumber\\
& \ge &\delta\big[s_n(\pi)-  \stn \1(0<|
 f^{(1)}(X_t)- f^{(2)}(X_t)|\le \delta)\big]\,.
\end{eqnarray}
Taking expectations on both sides of~(\ref{EQ:ISR_pr1}), we obtain that $R_n(\pi)\ge \delta\big[S_n(\pi)- n\delta^\alpha\big]$,
where we used the margin condition. The proof follows by choosing $\delta=(S_n(\pi)/cn)^{1/\alpha}$ for $c\ge2$ large enough to ensure
that $\delta <\delta_0$ \epr
Using Lemma~\ref{LEM:ISR}, we obtain the following corollary of
Theorem~\ref{TH:main}

\begin{cor}\label{COR}
Fix $\beta \in (0,1]$, $L>0$ and $\alpha \in (0,1]$. Under the conditions of Theorem~\ref{TH:main}, the UCBogram policy $\hpi$ with
$M=\lfloor(n/\log n)^{1/(2\beta+d)}\rfloor$ has an inferior sampling rate at time $n$
bounded by
$$
S_n(\hat \pi)\le Cn\Big(\frac{n}{\log
  n}\Big)^{-\frac{\beta\alpha}{2\beta+d}}\,.
$$
where $C>0$ is a positive constant.
\end{cor}

\section{Lower bound}
\label{SEC:LB}
While the UCBogram is a very simple policy, it still provides good
insights as to how to construct a lower bound on the regret for
incurred by any admissible policy.  Indeed, the main result of
this section demonstrates the polynomial rate of the upper bounds in Theorem~\ref{TH:main} and Corollary~\ref{COR} is optimal in a
minimax sense, for a large class of conditional reward
distributions. Define the Kullback-Leibler (KL) divergence between $P$ and $Q$,
where $P$ and $Q$ are two probability distributions by
$$
\cK(P,Q)=\left\{\begin{array}{ll}
\int\log\left(\frac{\ud P}{\ud Q}\right)\ud P & {\rm if} \ P \ll Q\,,\\
\infty &{\rm otherwise.}
\end{array}\right.
$$
Denote by $P\bi_{f(X)}$ the conditional distribution of
$Y\bi$ given $X$ for any $i=1,2$ and assume that there exists $\kappa^2>0$ such that for any $\theta, \theta' \in \Theta$ the KL divergence between $P_\theta\bi$ and $P_{\theta'}\bi$ satisfies
\begin{equation}
\label{eq:KL}
\cK(P_\theta\bi, P_{\theta'}\bi)\le \frac{1}{\kappa^2} (\theta-\theta')^2\,.
\end{equation}
Assumption~\eqref{eq:KL} is similar to Assumption~(B) employed in
\citet[Section 2.5]{Tsy09} but does not require absolute continuity
with respect to the Lebesgue measure. A direct consequence of the
following lemma is that Assumption~\eqref{eq:KL} is satisfied when $P_\theta$ is
a Bernoulli distribution with parameter $\theta \in (0,1)$.

\begin{lem}
For any  $a \in [0,1]$ and $b \in (0,1)$ let $P_a$ and $P_b$ denote two Bernoulli
distributions with parameters $a$ and $b$ respectively. Then
$$
\cK(P_a, P_b)\le \frac{(a-b)^2}{b(1-b)}\,.
$$
In particular, if  $b_0 \in [0,1/2)$,
Assumption~\eqref{eq:KL} is satisfied with $\kappa^2=1/4-b_0^2$, for any $a \in [0,1], b \in [1/2-b_0, 1/2+b_0]$.
\end{lem}
{\bf Proof.} From the definition of the KL divergence, we have
\begin{equation*}
\begin{split}
\cK(P_a, P_b)&=a\log\Big(\frac{a}{b} \Big)+
(1-a)\log\Big(\frac{1-a}{1-b} \Big)\le a\Big(\frac{a-b}{b} \Big)-(1-a)\Big(\frac{a-b}{1-b} \Big)= \frac{(a-b)^2}{b(1-b)}\\
\end{split}
\end{equation*}
where in the second line we used the inequality $\log(1+u)\le u$.

\epr

\begin{TH1}
\label{th:minimax} Fix $\alpha, \beta, L>0$ such that
$\alpha\beta<1$ and let $\cX=[0,1]^d$. Assume that the covariates
$X_t$ are uniformly distributed on the unit hypercube $\cX$ and
that there exists $\tau\in (0,1/2)$ such that $\{P\bi_{\theta}\,,
\ \theta \in [1/2-\tau,1/2+\tau]\}$ satisfies
equation~(\ref{eq:KL}) for $i=1,2$. Then, there exists a pair of
reward functions $f\bi, i=1, 2$ that satisfy both the smoothness
condition with parameters $(\beta, L)$ and the margin condition
with parameter $\alpha$,  such that for \emph{any}
non-anticipating policy $\pi$ the regret is bounded as follows
\begin{equation}
\label{EQ:LBreg}
R_n(\pi)\ge Cn^{1-\frac{\beta(\alpha+1)}{2\beta+d}}\,,
\end{equation}
and the inferior sampling rate is bounded as follows
\begin{equation}
\label{EQ:LBisr}
S_n(\pi)\ge Cn^{1-\frac{\beta\alpha}{2\beta+d}}\,,
\end{equation}for some positive constant $C$.
\end{TH1}
{\bf Proof.}  To
simplify the arguments below, it will be useful to denote arm $2$
by $-1$. Finally, with slight abuse of notation, we use
 $S_n(\pi, f^{(1)}, f^{(-1)})$
to denote the inferior sampling rate at time $n$ that is defined
in~\eqref{e-isr}, making the dependence on the mean reward
functions explicit.

In view of Lemma~\ref{LEM:ISR}, it is sufficient to
prove~\eqref{EQ:LBisr}. To do so we reduce our problem to a
hypothesis testing problems; an approach this is quite standard in
the nonparametric literature, cf.  \cite[Chapter~2]{Tsy09}. For
any policy $\pi$, and any $t=1, \ldots, n$, denote by
$\p_{\pi,f}^t$ the joint distribution of the collection of pairs
$$(X_1, Y_1^{(\pi_1(X_1))}), \ldots, (X_t, Y_t^{(\pi_t(X_t))})$$
where $\E[Y^{(1)}|X]=f(X)$ and $\E[Y^{(-1)}|X]=1/2$. Let
$\E_{\pi,f}^t$ denote the corresponding expectation. It follows
that the oracle policy $\spi_f$ is given by $\spi_f(x)=\sgn[f(x)]$
with the convention that $\sgn(0)=1$. Fix $\delta_0 \in (0,1)$ as
in the definition of the margin condition. We now construct a
class $\cC$ of functions $f:\cX \to [0,1]$ such that $f$
satisfies~\eqref{EQ:holder} and
$$
P_X\big[ \,0< |f(X)-1/2 |\le \delta\big]\le
C_\delta \delta^\alpha\,, \quad \forall\,   \delta \in [0, \delta_0]\,,
$$
As a result, the machine characterized by the expected rewards
$f^{(1)}=f$ and $f^{(-1)}=1/2$ satisfies both the smoothness and
the margin conditions. Moreover, we construct $\cC$ in such a way
that for any policy $\pi$
\begin{equation}
\label{EQ:LBproof}
\sup_{f \in \cC}S_n(\pi, f, 1/2)\ge Cn\Big(\frac{n}{\log n}\Big)^{-\frac{\beta\alpha}{2\beta+d}}\,.
\end{equation}
for some positive constant $C$.
%
%
Consider the
 regular grid
$\cQ=\{q_1, \ldots, q_{M^d}\}$, where $q_k$ denotes the center of bin
$B_k$, $k=1, \ldots, M^d$, for some $M\ge 1$ to be defined. Define $C_\phi=\min(L, \tau, 1/4)$  and let $\phi_\beta:  \R^d \to \R_+ $  be a smooth function defined as follows:
$$
\phi_\beta(x)=\left\{
\begin{array}{ll}
 (1-\|x\|_\infty)^\beta & {\rm if} \ 0\le \|x\|_\infty \le 1,\\
0 & {\rm if} \ \|x\|_\infty>1\,.\\
\end{array}
\right.
$$
Clearly, we have
$|C_\phi\phi_\beta(x)-C_\phi\phi_\beta(x')| \le L \|x-x'\|_\infty^\beta\le L \|x-x'\|^\beta$ for any $x,
x'\in \R^d$.

Define the integer $m = \lceil \mu M^{d-\alpha \beta}\rceil$,
i.e., the smallest integer that is larger than or equal to $\mu
M^{d-\alpha \beta}$, where $\mu \in (0,1)$ is chosen small enough
to ensure that $m\le M^d$. Define $\Omega_m=\{-1,1\}^m$ and for
any $\omega\in \Omega_m$, define the function $f_\omega$ on
$[0,1]^d$ by
$$
f_{\omega}(x)=1/2+\sum_{j=1}^{m} \omega_j \varphi_j(x)\,,
$$
where $\varphi_j(x)=M^{-\beta}C_\phi \phi (M[x-q_j])\1(x \in
  B_j)$. Notice in particular that $f_\omega(x)=1/2$ if and only if $x \in
\cX\setminus \bigcup_{j=1}^m B_j$ up to a set of zero Lebesgue
measure. We are now in position to define the family $\cC$ as
$$
\cC = \left\{ f_{\omega}\;:\; \omega \in \Omega_m\right\}\,.
$$
Note first that any function $f_\omega \in \cC$ satisfies the
smoothness condition~\eqref{EQ:holder}. We now check that the
margin condition is satisfied with parameter $\alpha$. For any
$\omega \in \Omega_m$, we have
\begin{equation*}
\label{EQ6:pr:thlb1}
\begin{split}
P_X(0<|f_{\omega}(X)-1/2| \le C_\phi \delta )&= \sum_{j=1}^mP_X(0<|f_{\omega}(X)-1/2| \le
C_\phi\delta, X \in B_j)\\
&= m P_X(0< \phi (M[X-q_1])\le \delta M^{\beta}, X \in B_1)\,\\
&= m \int_{B_1}\hspace{-.3em}\1(\phi (Mx) \le \delta M^{\beta})\ud x \,\\
&= mM^{-d}\int_{[0,1]^d}\hspace{-.8em} \1(\phi (x) \le \delta
  M^{\beta})\ud x\,,
\end{split}
\end{equation*}
where in the third equality, we used the fact that $P_X$ denotes the
uniform distribution on $[0,1]^d$.
Now, since $\phi$ is non negative and uniformly bounded by 1, we have on the one hand that for $\delta M^\beta>1$,
$$
\int_{[0,1]^d}\hspace{-.8em} \1(\phi (x) \le \delta
  M^{\beta})\ud x =1\,.
$$
On the other hand, when $\delta M^\beta\le 1$, we find
$$
\int_{[0,1]^d}\hspace{-.8em} \1(\phi (x) \le \delta
  M^{\beta})\ud x=1-\int_{[0,1]^d}\hspace{-.8em}\1(\|x\|_\infty \le 1-M\delta^{1/\beta}
  )\ud x=1-\left(1-M\delta^{1/\beta}
  \right)^d \le dM\delta^{1/\beta}\,.
$$
It yields
\begin{equation*}
\label{EQ6:pr:thlb2}
\begin{split}
P_X(0<|f_{\omega}(X)-1/2| \le C_\phi \delta )&\le 
mM^{-d}\1(\delta M^{\beta} > 1)+
mdM^{1-d}\delta^{1/\beta} \1(\delta M^{\beta}  \le 1)\big) \\
&\le M^{-\alpha\beta} \1( M^{-\alpha
    \beta}<\delta^\alpha) +dM^{1-\alpha\beta}\delta^{1/\beta}
\1(M\le \delta^{-1/\beta}) \\
    & \le (1+d) \delta^{\alpha}\,,\\
\end{split}
\end{equation*}
where  we used the fact that $1-\alpha \beta \ge 0$ to bound the
second term in the last inequality. Thus, the margin condition is
satisfied for any $\delta_0$ and with $C_\delta=(1+d)/C_\phi^\alpha$.

We now prove~\eqref{EQ:LBproof} by observing that if we denote $\omega=(\omega_1, \ldots, \omega_m) \in \Omega_m$, we have
\begin{eqnarray}
\sup_{f \in \cC}S_n(\pi, f^{(1)}, 1/2) & = & \sup_{\omega \in \Omega_m} \sum_{t=1}^n \E_{\pi,f_\omega}^{t-1} P_X\left[\pi_t(X_t)\neq \sgn(f_\omega(X_t))\right]\nonumber \\
&=& \sup_{\omega \in \Omega_m} \sum_{j=1}^m\sum_{t=1}^n \E_{\pi,f_\omega}^{t-1}P_X\left[\pi_t(X_t)\neq \omega_j , X_t \in B_j\right]\nonumber\\
&\ge& \frac{1}{2^m}  \sum_{j=1}^m\sum_{t=1}^n \sum_{\omega \in \Omega_m}\E_{\pi,f_\omega}^{t-1}P_X\left[\pi_t(X_t)\neq \omega_j , X_t \in B_j\right] \label{EQ:prLBsn}
\end{eqnarray}
Observe now that for any $j=1, \ldots, m$, the sum $\sum_{\omega
\in \Omega} [\cdots]$ in the previous display can be decomposed as
$$
Q_j^t=\sum_{\omega_{[-j]} \in \Omega_{m-1}}\sum_{i \in \{-1,1\}}\E_{\pi,f_{\omega_{[-j]}^i}}^{t-1}P_X\left[\pi_t(X_t)\neq i , X_t \in B_j\right]\,,
$$
where $\omega_{[-j]}=(\omega_1, \ldots, \omega_{j-1},
\omega_{j+1}, \ldots, \omega_m)$ and $\omega_{[-j]}^i=(\omega_1,
\ldots, \omega_{j-1}, i, \omega_{j+1}, \ldots, \omega_m)$ for
$i=-1,1$. Using Theorem~2.2$(iii)$ of \citet{Tsy09}, and denoting
by $P_X^j(\cdot)$ the conditional distribution $P_X(\cdot|X \in
B_j)$,  we get
\begin{eqnarray}
\label{EQ:arrayproof1}
\sum_{i \in \{-1,1\}}\E_{\pi,f_{\omega_{[-j]}^i}}^{t-1}P_X\left[\pi_t(X_t)\neq i , X_t \in B_j\right]&=&
\frac{1}{M^d}\sum_{i \in \{-1,1\}}\E_{\pi,f_{\omega_{[-j]}^i}}^{t-1}P_X^j\left[\pi_t(X_t)\neq i\right] \nonumber\\
&\ge& \frac{1}{4M^d}\exp\left[-\cK\big( \p_{\pi,f_{\omega_{[-j]}^{-1}}}^{t-1} \times P_X^j , \p_{\pi,f_{\omega_{[-j]}^1}}^{t-1}  \times P_X^j\big) \right] \nonumber\\
&=&\frac{1}{4M^d}\exp\left[-\cK\big( \p_{\pi,f_{\omega_{[-j]}^{-1}}}^{t-1}, \p_{\pi,f_{\omega_{[-j]}^1}}^{t-1} \big) \right]
\end{eqnarray}
For any $t=2, \ldots, n$, let $\cF_t$ denote the $\sigma$-algebra
generated by the information available at time $t$ immediately
{\it after} observing $X_t$, i.e., $\cF_t=\sigma\big(X_t,(X_s,
Y_s^{(\pi_s(X_s))}), s=1, \dots, t-1)\big)$. Define the
conditional distribution $\p_{\pi,f}^{\cdot|\cF_t}$ of the random
couple $(X_t, Y_t^{(\pi_t(X_t))})$, conditioned on $\cF_t$. Denote
also by $E_{X_t}$ the expectation with respect to the marginal
distribution of $X_t$.
 Applying the chain rule for KL divergence, we find that for any $t=1, \ldots, n$ and any $f, g : \cX \to [0,1]$, we have
\begin{eqnarray*}
\cK\big(\p_{\pi,f}^t, \p_{\pi,g}^t\big)&=&\cK\big(\p_{\pi,f}^{t-1}, \p_{\pi,g}^{t-1}\big)+\E_{\pi,f}^{t-1}E_{X_t}\left[\cK\big(\p_{\pi,f}^{\cdot|\cF_{t}}, \p_{\pi,g}^{\cdot|\cF_{t}}\big) \right]\\
&=&\cK\big(\p_{\pi,f}^{t-1}, \p_{\pi,g}^{t-1}\big)+\E_{\pi,f}^{t-1}E_{X_t}\left[\cK\big(\p_{\pi,f}^{Y_t^{(\pi_t(X_t))}|\cF_{t}}, \p_{\pi,g}^{Y_t^{(\pi_t(X_t))}|\cF_{t}}\big) \right]\,,
\end{eqnarray*}
where  $\p_{\pi,f}^{Y_t^{(\pi_t(X_t))}|\cF_{t}}$ denotes the
conditional distribution of $Y_t^{(\pi_t(X_t))}$ given $\cF_{t}$.
Since, for any $f \in \cC$, we have that
$\E[Y_t^{(\pi_t(X_t))}|\cF_{t}]=f^{(\pi_t(X_t))}(X_t) \in
[1/2-\tau, 1/2+\tau]$, we can apply~\eqref{eq:KL} to derive  the
following upper bound:
\begin{eqnarray*}
\cK\big( \p_{\pi,f_{\omega_{[-j]}^{-1}}}^{Y_t^{(\pi_t(X_t))}|\cF_{t}},\p_{\pi,f_{\omega_{[-j]}^{1}}}^{Y_t^{(\pi_t(X_t))}|\cF_{t}}\big) &\le& \frac{1}{\kappa^2} \left(f_{\omega_{[-j]}^{1}}(X_t)-f_{\omega_{[-j]}^{-1}}(X_t)\right)^2\1\left(\pi_t(X_t)=1\right) \\
&\le& \frac{4}{\kappa^2} C_\phi^2 M^{-2\beta}\1\left(\pi_t(X_t)=1, X_t \in B_j\right)  \\
&\le& \frac{M^{-2\beta}}{4\kappa^2}\1\left(\pi_t(X_t)=1, X_t \in B_j\right) \,.
\end{eqnarray*}
By induction, the last two displays yield that for any $t=1,
\ldots, n$,
\begin{equation}
\label{EQ:KL_UB}
\cK(\p_{\pi,f_{\omega_{[-j]}^1}}^{t-1} , \p_{\pi,f_{\omega_{[-j]}^{-1}}}^{t-1}) \le\frac{M^{-2\beta}}{4\kappa^2} \textsf{N}_{j, \pi}\,,
\end{equation}
where
$$
\textsf{N}_{j, \pi}=\E_{\pi,f_{\omega_{[-j]}^{-1}}}^{n-1}E_{X}\left[\sum_{t=1}^n \1\left(\pi_t(X)=1, X \in B_j\right)\right]\,,
$$
denotes the expected number of times $t$ between time 1 and time
$n$  that  $X_t \in B_j$ and $\pi_t(X_t)=1$.
Combining~\eqref{EQ:arrayproof1} and~\eqref{EQ:KL_UB}, we get
\begin{equation}
\label{EQ:Qj1}
Q_j^t\ge \frac{2^{m-1}}{4M^d}\exp\left(-\frac{M^{-2\beta}}{4\kappa^2 }\textsf{N}_{j,\pi}\right)\,.
\end{equation}
On the other hand, from the definition of $Q_j^t$, we clearly have
\begin{equation}
\label{EQ:Qj2}
\sum_{t=1}^n  Q_j^t \ge 2^{m-1}\textsf{N}_{j,\pi}\,.
\end{equation}
Plugging the lower bounds~\eqref{EQ:Qj1}  and~\eqref{EQ:Qj2} into~\eqref{EQ:prLBsn} yields
\begin{eqnarray*}
\sup_{f \in \cC}S_n(\pi, f^{(1)}, 1/2)& \ge & \frac{2^{m-1}}{2^m}\sum_{j=1}^m \max\left\{\frac{n}{4M^d} \exp\left(-\frac{M^{-2\beta}}{4\kappa^2 }\textsf{N}_{j,\pi}\right) , \textsf{N}_{j,\pi}\right\}\\
& \ge & \frac{1}{4} \sum_{j=1}^m \left\{\frac{n}{4M^d} \exp\left(-\frac{M^{-2\beta}}{4\kappa^2 }\textsf{N}_{j,\pi}\right) + \textsf{N}_{j,\pi}\right\}\\
& \ge &\frac{m}{4}  \inf_{z\ge 0}\left\{\frac{n}{4M^d} \exp\left(-\frac{M^{-2\beta}}{4\kappa^2 }z\right) + z\right\}
\end{eqnarray*}
Notice now that
$$
z^*=\argmin_{z\ge 0}\left\{\frac{n}{4M^d} \exp\left(-\frac{M^{-2\beta}}{4\kappa^2 }z\right) + z\right\}
$$
is strictly positive if and only if $n > 16\kappa^2 M^{2\beta+d}$, in which case
$$
z^*=4\kappa^2 M^{2\beta}\log\left( \frac{n}{16\kappa^2M^{2\beta+d}}\right)\,.
$$
Taking
$$
M=\left\lceil \left(\frac{ n}{16e\kappa^2}\right)^{\frac{1}{2\beta+d}} \right\rceil
$$
gives $z^*=c^*n^{\frac{2\beta}{2\beta+d}} $ for some positive constant $c^*$, so that
$$
\sup_{f \in \cC}S_n(\pi, f^{(1)}, 1/2)\ge Cmz^* \ge
Cn^{1-\frac{\alpha\beta}{2\beta+d}}.
$$
This completes the proof.  \epr

Notice that the rates obtained in Theorem~\ref{th:minimax}, can be obtained in the full information case, where the operator observes the whole i.i.d sequence $(X_i, Y^{(1)}_i, Y^{(2)}_i), i=1, \ldots, n$, even before the first round. Indeed, such bounds have been obtained by~\citet{AudTsy07} in the classification setup, i.e., when the rewards are Bernoulli random variables. However, we state a different technique, tailored for bandit policies in a partial information setup. While the final result is the same, we believe that it sheds light on the technicalities encountered in proving such a lower bound.

\bigskip

\begin{flushright}

\footnotesize

{\sc Philippe Rigollet}\\

{\sc Department of Operations Research} \\

{\sc and Financial Engineering}\\

{\sc Princeton University}\\

{\sc Princeton, NJ 08544, USA}\\

\texttt{  rigollet@princeton.edu}

\bigskip

{\sc Assaf Zeevi}\\

{\sc Graduate School of Business}\\
{\sc Columbia University}\\
{\sc New York, NY 10027}\\
\texttt{ assaf@gsb.columbia.edu}\\

\end{flushright}


\begin{thebibliography}{19}
\expandafter\ifx\csname natexlab\endcsname\relax\def\natexlab#1{#1}\fi
\expandafter\ifx\csname url\endcsname\relax
  \def\url#1{\texttt{#1}}\fi
\expandafter\ifx\csname urlprefix\endcsname\relax\def\urlprefix{URL }\fi
\providecommand{\eprint}[2][]{\url{#2}}

\bibitem[{Audibert and Tsybakov(2007)}]{AudTsy07}
\textsc{Audibert, J.-Y.} and \textsc{Tsybakov, A.~B.} (2007).
\newblock Fast learning rates for plug-in classifiers.
\newblock \textit{Ann. Statist.}, \textbf{35} 608--633.

\bibitem[{Auer et~al.(2002)Auer, Cesa-Bianchi and Fischer}]{AueCesFis02}
\textsc{Auer, P.}, \textsc{Cesa-Bianchi, N.} and \textsc{Fischer, P.} (2002).
\newblock Finite-time analysis of the multiarmed bandit problem.
\newblock \textit{Mach. Learn.}, \textbf{47} 235--256.

\bibitem[{Boucheron et~al.(2005)Boucheron, Bousquet and Lugosi}]{BouBouLug05}
\textsc{Boucheron, S.}, \textsc{Bousquet, O.} and \textsc{Lugosi, G.} (2005).
\newblock Theory of classification: a survey of some recent advances.
\newblock \textit{ESAIM Probab. Stat.}, \textbf{9} 323--375 (electronic).

\bibitem[{Cesa-Bianchi and Lugosi(2006)}]{CesLug06}
\textsc{Cesa-Bianchi, N.} and \textsc{Lugosi, G.} (2006).
\newblock \textit{Prediction, learning, and games}.
\newblock Cambridge University Press, Cambridge.

\bibitem[{Goldenshluger and Zeevi(2009)}]{GolZee09}
\textsc{Goldenshluger, A.} and \textsc{Zeevi, A.} (2009).
\newblock Woodroofe's one-armed bandit problem revisited.
\newblock \textit{Ann. Appl. Probab.}, \textbf{19} 1603--1633.

\bibitem[{Hazan and Megiddo(2007)}]{HazMeg07}
\textsc{Hazan, E.} and \textsc{Megiddo, N.} (2007).
\newblock Online learning with prior knowledge.
\newblock In \textit{Learning theory}, vol. 4539 of \textit{Lecture Notes in
  Comput. Sci.} Springer, Berlin, 499--513.

\bibitem[{Kakade et~al.(2008)Kakade, Shalev-Shwartz and Tewari}]{KakShaTew08}
\textsc{Kakade, S.}, \textsc{Shalev-Shwartz, S.} and \textsc{Tewari, A.}
  (2008).
\newblock Efficient bandit algorithms for online multiclass prediction.
\newblock In \textit{Proceedings of the 25th Annual International Conference on
  Machine Learning (ICML 2008)} (A.~McCallum and S.~Roweis, eds.). Omnipress,
  440--447.

\bibitem[{Lai and Robbins(1985)}]{LaiRob85}
\textsc{Lai, T.~L.} and \textsc{Robbins, H.} (1985).
\newblock Asymptotically efficient adaptive allocation rules.
\newblock \textit{Adv. in Appl. Math.}, \textbf{6} 4--22.

\bibitem[{Langford and Zhang(2008)}]{LanZha08}
\textsc{Langford, J.} and \textsc{Zhang, T.} (2008).
\newblock The epoch-greedy algorithm for multi-armed bandits with side
  information.
\newblock In \textit{Advances in Neural Information Processing Systems 20}
  (J.~Platt, D.~Koller, Y.~Singer and S.~Roweis, eds.). MIT Press, Cambridge,
  MA, 817--824.

\bibitem[{Lu et~al.(2009)Lu, P{\'a}l and P{\'a}l}]{LuPalPal09}
\textsc{Lu, T.}, \textsc{P{\'a}l, D.} and \textsc{P{\'a}l, M.} (2009).
\newblock Showing relevant ads via context multi-armed bandits.
\newblock Tech. rep.

\bibitem[{Mammen and Tsybakov(1999)}]{MamTsy99}
\textsc{Mammen, E.} and \textsc{Tsybakov, A.~B.} (1999).
\newblock Smooth discrimination analysis.
\newblock \textit{Ann. Statist.}, \textbf{27} 1808--1829.

\bibitem[{Rigollet and Vert(2009)}]{RigVer09}
\textsc{Rigollet, P.} and \textsc{Vert, R.} (2009).
\newblock Fast rates for plug-in estimators of density level sets.
\newblock \textit{Bernoulli}, \textbf{15} 1154--1178.

\bibitem[{Robbins(1952)}]{Rob52}
\textsc{Robbins, H.} (1952).
\newblock Some aspects of the sequential design of experiments.
\newblock \textit{Bull. Amer. Math. Soc.}, \textbf{58} 527--535.

\bibitem[{Slivkins(2009)}]{Sli09}
\textsc{Slivkins, A.} (2009).
\newblock Contextual bandits with similarity information.
\newblock \textit{Arxiv preprint arXiv:0907.3986}.

\bibitem[{Tsybakov(2004)}]{Tsy04}
\textsc{Tsybakov, A.~B.} (2004).
\newblock Optimal aggregation of classifiers in statistical learning.
\newblock \textit{Ann. Statist.}, \textbf{32} 135--166.

\bibitem[{Tsybakov(2009)}]{Tsy09}
\textsc{Tsybakov, A.~B.} (2009).
\newblock \textit{Introduction to Nonparametric Estimation}.
\newblock Springer Publishing Company, Incorporated.

\bibitem[{Wang et~al.(2005)Wang, Kulkarni and Poor}]{WanKulPoo05}
\textsc{Wang, C.-C.}, \textsc{Kulkarni, S.} and \textsc{Poor, H.} (2005).
\newblock Bandit problems with side observations.
\newblock \textit{Automatic Control, IEEE Transactions on}, \textbf{50}
  338--355.

\bibitem[{Woodroofe(1979)}]{Woo79}
\textsc{Woodroofe, M.} (1979).
\newblock A one-armed bandit problem with a concomitant variable.
\newblock \textit{J. Amer. Statist. Assoc.}, \textbf{74} 799--806.

\bibitem[{Yang and Zhu(2002)}]{YanZhu02}
\textsc{Yang, Y.} and \textsc{Zhu, D.} (2002).
\newblock Randomized allocation with nonparametric estimation for a multi-armed
  bandit problem with covariates.
\newblock \textit{Ann. Statist.}, \textbf{30} 100--121.

\end{thebibliography}
\end{document}